\documentclass[11pt]{article}
\usepackage{graphicx}
\usepackage{amsmath,amsthm,amssymb}
\usepackage{xypic}
\usepackage{mathabx}
\usepackage{bbm}
\usepackage{color}
\newtheorem{theorem}{Theorem}

\newtheorem{example}[theorem]{Example}

\newtheorem{remark}[theorem]{Remark}
\newtheorem{corollary}[theorem]{Corollary}

\begin{document}

\title{Euler Characteristics of Finite Homotopy Colimits}
\author{John D. Berman}
\maketitle

\begin{abstract}
We provide a brief calculation of the Euler characteristic of a finite homotopy colimit of finite cell complexes, which depends only on the Euler characteristics of each space and resembles Mobius inversion. Versions of the result are known when the colimit is indexed by a finite category, but the behavior is more uniform when we index by finite quasicategories instead. The formula simultaneously generalizes the additive formula for Euler characteristic of a homotopy pushout and the multiplicative formula for Euler characteristic of a fiber bundle.
\end{abstract}

\noindent Given a fibration of finite cell complexes $E\rightarrow B$ with $B$ connected and whose fiber $F$ is also finite, the Euler characteristic satisfies the formula $\chi(E)=\chi(F)\chi(B)$. This is a classical fact with a number of proofs, notably via the Serre Spectral Sequence. In this note, we present a short and explicit proof of a more general result. That is, we describe a Mobius inversion formula for computing Euler characteristics of finite homotopy colimits.

This formula is hard to trace back, but versions of it go back decades. Recently, it has been studied by Leinster \cite{Leinster} and Fiore-L\"uck-Sauer \cite{FLS1}\cite{FLS2}, attached to names like the \emph{Euler characteristic} or \emph{Mobius function} of categories. However, it does not appear to be written down at the generality of $\infty$-categories. We want to correct this, first because it is such a (comparatively) simple application of the quasicategorical technology that it ought to be better known; second because it is more general than the 1-categorical statement. For example, the fibration formula is a special case.

\begin{remark}
By a `finite homotopy colimit', we mean a homotopy colimit indexed by a finite $\infty$-category. Finite $\infty$-categories are very different from finite categories (a fact which may be surprising when first encountered). For example, if $BG$ denotes the category with one object associated to a group $G$, $B\mathbb{Z}/2\cong\mathbb{R}P^\infty$ is finite as a category but not as an $\infty$-category, while $B\mathbb{Z}\cong S^1$ is finite as an $\infty$-category but not as a category.

This is one reason why some treatments of the Mobius inversion formula require restrictive assumptions on the indexing category.
\end{remark}

\noindent We will state our main result first as a model-independent (but nonconstructive) fact about Mobius functions on $\infty$-categories. Then we will describe how to compute Mobius functions, given a finite cellular model for the $\infty$-category. This is a categorification of the usual alternating sum formula for the Euler characteristic.

\begin{theorem}
\label{Thm1}
If $\mathcal{K}$ is a finite $\infty$-category, let $\bar{\mathcal{K}}$ denote the set of equivalence classes of objects . There is a \emph{Mobius function} $\mu:\bar{\mathcal{K}}\rightarrow\mathbb{Z}$ with the property: For any functor $\mathcal{K}\xrightarrow{f}\text{Top}$ which takes values in finite cell complexes, $$\chi(\text{colim}\,f)=\sum_{[x]\in\bar{\mathcal{K}}}\mu[x]\chi(f[x]).$$
\end{theorem}

\noindent Moreover, given a model for $\mathcal{K}$ as a finite simplicial set (in the Joyal model structure), the Mobius function can be directly computed as follows:

\begin{theorem}
\label{Thm2}
A simplicial set is called \emph{finite} if it has finitely many nondegenerate simplices. If $K$ is a finite simplicial set, and $x\in K_0$ is an object (0-simplex), the Mobius function is $\mu(x)=\sum(-1)^n|K_n^x|$. Here $|K_n^x|$ denotes the number of nondegenerate $n$-simplices with initial vertex $x$.

If $K\xrightarrow{f}\text{Top}$ is a functor that lands in finite cell complexes, $$\chi(\text{colim}\,f)=\sum_{x\in K_0}\mu(x)\chi(f(x)).$$
\end{theorem}

\begin{remark}
Note that $\mu$ from Theorem \ref{Thm1} is applied to equivalence classes $[x]$ of objects, while $\mu$ from Theorem \ref{Thm2} is applied to objects themselves. If $K$ is a finite simplicial set modeling an $\infty$-category, and $[x]$ an equivalence class of objects (0-simplices) the two are related via $$\mu[x]=\sum_{x\in[x]}\mu(x).$$
\end{remark}

\noindent Many classical properties of Euler characteristics are special cases:

\begin{example}[$\chi$ of a simplicial complex]
\label{Ex1}
For a finite simplicial set (or finite $\infty$-category) $K$, the classifying space $|K|$ is the colimit of the constant diagram at a point $K\xrightarrow{\ast}\text{Top}$ (\cite{HTT} 4.4.4.9). Therefore, the Euler characteristic of its classifying space is given by $$\chi|K|=\sum_{x\in K_0}\mu(x)=\sum_{n\geq 0}(-1)^n|K_n|.$$ We recover the classical fact that the Euler characteristic of a simplicial complex is the alternating sum of the numbers of $n$-simplices.
\end{example}

\begin{example}[$\chi$ of a homotopy pushout]
\label{Ex2}
Given a homotopy pushout $$\xymatrix{
A\ar[r]\ar[d] &B\ar[d] \\
C\ar[r] &D
}$$ of finite cell complexes, then $\chi(D)=\chi(B)+\chi(C)-\chi(A)$.
\end{example}

\begin{example}[$\chi$ of a fiber bundle]
Regard a finite cell complex $B$ as an $\infty$-groupoid, also finite by \cite{HTT} 1.2.14.2. If $B$ is connected, then the $\infty$-groupoid has a single object $b$ up to equivalence, and $\mu(b)=\chi(B)$ by Example \ref{Ex1}.

For any map $E\rightarrow B$ with fiber $F$, there is a functor $E_{-}:B\rightarrow\text{Top}$ which sends a point in $B$ to the homotopy fiber over that point, and $E_x\cong F$ for all $x\in B$ since $B$ is connected. Moreover, $\text{colim}(E_{-})\cong E$ by \cite{HTT} 3.3.4.6. Apply Theorem \ref{Thm1} with $\mu(b)=\chi(B)$:
\end{example}

\begin{corollary}
If $F\rightarrow E\rightarrow B$ is a fiber sequence of finite CW complexes with $B$ connected, then $\chi(E)=\chi(F)\chi(B)$.
\end{corollary}

\noindent We also offer a more elaborate example:

\begin{example}
Consider the diagram category $I$ of the form $$\xymatrix{
A\ar[r]\ar[d] &B\ar@<-.5ex>[r] \ar@<.5ex>[r] &C\\
D, &&
}$$ where the two composites $A\rightarrow C$ are equal. If $I\rightarrow\text{Top}$ is a diagram picking out finite cell complexes $A,B,C,D$, then the homotopy colimit has Euler characteristic $$\chi(\text{colim})=\chi(C)+\chi(D)-\chi(B)-2\chi(A).$$
\end{example}

\noindent We end by proving Theorem \ref{Thm2} (which implies Theorem \ref{Thm1}). The proof exactly follows Lurie's proof (for $\infty$-categories) that all finite colimits can be built out of pushouts and initial objects, \cite{HTT} 4.4.2.4. We take for granted that the Euler characteristic is 0 on the empty space and takes homotopy pushouts to sums (as in Example \ref{Ex2}).

Choose $K$ finite and $K\xrightarrow{f}\text{Top}$ which lands in finite cell complexes. We induct first on the dimension of $K$, then on the number of simplices of maximal dimension. When $K$ is 0-dimensional, the theorems hold by $\chi(X\amalg Y)=\chi(X)+\chi(Y)$.

Otherwise, fix a nondegenerate simplex $\Delta^n$ in $K$ of maximal dimension, with initial vertex $s$ and terminal vertex $t$. Decompose $K\cong K^\prime\amalg_{\partial\Delta^n}\Delta^n$, and let $\mu$ be the Mobius function on $K$, $\mu^\prime$ the Mobius function on $K^\prime$ (as defined in Theorem \ref{Thm2}). Note that $K$ has the same objects as $K^\prime$ since $n\geq 1$. Moreover, $\mu(x)=\mu^\prime(x)$, except at $x=s$, where $\mu(s)=\mu^\prime(s)+(-1)^n$.

Now we compute the colimit. Let $$X=\text{colim}(f|_{K^\prime})$$ $$Y=\text{colim}(f|_{\partial\Delta^n})$$ $$Z=\text{colim}(f|_{\Delta^n}).$$ Then $\text{colim}(f)=X\amalg_YZ$, so $\chi(\text{colim}\,f)=\chi(X)+\chi(Z)-\chi(Y)$. By the induction hypothesis, $\chi(Y)=\chi(f(t))-(-1)^n\chi(f(s))$ and $\chi(Z)=\chi(f(t))$, so $$\chi(\text{colim}\,f)=\chi(X)+(-1)^n\chi(f(s))=\sum_x\mu^\prime(x)\chi(f(x))+(-1)^n\chi(f(s))$$ Since $\mu(x)=\mu^\prime(x)$ and $\mu(s)=\mu^\prime(s)+(-1)^n$, this completes the proof.

\begin{remark}
Notice that we have proven something slightly more general. Suppose $\mathcal{C}$ is an $\infty$-category which admits finite colimits, $A$ is an abelian group, and $\chi:\mathcal{C}\rightarrow A$ is a function satisfying the three conditions:
\begin{itemize}
\item if $X\cong Y$, then $\chi(X)=\chi(Y)$;
\item if $0$ is the initial object, then $\chi(0)=0$;
\item if $$\xymatrix{
A\ar[r]\ar[d] &B\ar[d] \\
C\ar[r] &D
}$$ is a homotopy pushout square, then $\chi(A)+\chi(D)=\chi(B)+\chi(C)$.
\end{itemize}
Then for any finite $\infty$-category $K$ and functor $K\xrightarrow{f}\mathcal{C}$, we have $$\chi(\text{colim}\,f)=\sum_{[x]\in\bar{K}}\mu[x]\chi(f[x]).$$
\end{remark}

\end{document}